%% file: Harmonic_v3.tex
\documentclass[final,11pt]{article}

\usepackage{graphicx}
\usepackage{amsmath}
\usepackage{xcolor}
\usepackage{cite}
\usepackage{epstopdf} 
\usepackage{mathptmx}      

\newtheorem{theorem}{Theorem}
\newtheorem{lemma}{Lemma}
\newtheorem{remark}{Remark}

\newtheorem{example}{Example}[section]

\input{zmacros}

\begin{document}

\title{Analysis of harmonic average method for interface problems with discontinuous solutions and fluxes}

\author{Kejia Pan\thanks{School of Mathematics and Statistics, HNP-LAMA, Central South University, Changsha, Hunan, 410083, China, Email: kejiapan@csu.edu.cn }
\and
Hengrui Xu\thanks{School of Mathematics and Statistics, HNP-LAMA, Central South University, Changsha, Hunan, 410083, China }
\and Zhilin Li \thanks{Department of Mathematics, North Carolina State University, Raleigh, NC 27695-8205, USA, Email: zhilin@math.ncsu.edu}
}

\maketitle

\begin{abstract}
Harmonic average method has been widely utilized to deal with heterogeneous coefficients in solving differential equations. One remarkable advantage of the harmonic averaging method is that no derivative of the coefficient is needed. Furthermore, the coefficient matrix of the finite difference equations is an M-matrix which guarantees the stability of the algorithm.
It has been numerically observed but not theoretically proved that the method produces second order pointwise accuracy when the solution and flux are continuous even if the coefficient has finite discontinuities for which  the method is inconsistent ($O(1)$ in the local truncation errors). It has been believed that there are some fortunate error cancellations. The  harmonic average method does not converge when the solution or the flux has finite discontinuities.
In this paper, not only we rigorously prove the second order convergence of the harmonic averaging method for one-dimensional interface problem when the coefficient has a finite discontinuities and the solution and the flux are continuous, but also proposed an {\em  improved harmonic average method} that is also second order accurate (in the $L^{\infty}$ norm), which allows discontinuous solutions and fluxes  along with the discontinuous coefficients. The key in the convergence proof is the construction of the  Green's function. The proof shows how the error cancellations occur in a subtle way. Numerical experiments in both 1D and 2D confirmed the theoretical proof of the improved harmonic average method. 

\end{abstract}

{\bf keywords:}
 Harmonic average, improved harmonic average method, variable discontinuous coefficient, non-homogeneous jump conditions, Green function,  discrete maximum principle, convergence analysis.


\section{Introduction}

Harmonic average, sometimes also called harmonic mean, has been applied to solve various differential equations when the material parameters have large variations, even with finite jump discontinuities, which are called interface problems.
Applications include unsaturated flow in layered soils \cite{MR3771185},  non-linear heat conduction problems in \cite{osti_928087}. Advantages of the harmonic average approach include the simplicity in implementation and 
preservation of some physical properties. There are limited references on study of the harmonic average method. A few discussions can be found in \cite{szy-helm11,MR3771185,MR4762021,osti_928087}.  It is observed that the harmonic average method works well for one-dimensional problems when the solution and the flux are continuous even if the coefficient have a finite number of jump discontinuities. However, to our best knowledge, there is no rigorous proof that can be found in the literature. In this paper, we have provided rigorous proof why and when the harmonic average method works; and more important,  developed an improved harmonic average method that can work for interface problems with discontinuous solution and/or fluxes with second order accuracy in the $L^{\infty}$ norm. The new  improved harmonic average method does not need the derivative of the coefficient, an obvious advantage over some existing second order accurate methods.

We consider  the following interface problem
\eqml{eq1d}
 \dsp (\beta(x) u_x)_x - \sigma(x) u = f(x) + v \delta(x-\alf) + w \delta'(x-\alf) , \quad  0< x < 1,
\enml
with specified boundary conditions of $u(x)$ at $x=0$ and $x=1$.  We assume that $\beta(x)\ge \beta_0>0$, $\sigma(x)\ge 0$,  and $\sigma(x), f(x)\in C[0,1]$, but allow $\beta(x) \in C(0,\alf)  \cup C(\alf,1)$, which means that $\beta(x)$ can have a finite jump discontinuity at $\alf$. In \cite{REU10}, the authors showed that the problem is equivalent to the following interface problem,
\eqml{eq1dB}
& \dsp (\beta(x) u_x)_x - \sigma(x) u = f(x), \qquad \{0 < x < \alf\} \cup \{  \alf < x <1\}, \\ \eqsp
& \dsp [u]_{\alf} = \frac{2 w}{\beta^- + \beta^+}, \qquad [\beta u_x] _{\alf} = v, \qquad 0 < \alf < 1,
\enml
where $\beta^-= \lim_{x\goto \alf-} \beta(x)$ and $\beta^+= \lim_{x\goto \alf+} \beta(x)$. It is easier to use this formulation to prove the existence and uniqueness of the solution to the boundary value problem.

\vthin

{\bf Harmonic average finite difference method.} Given a uniform grid $0=x_0<x_1<\cdots<x_{N}=1$ with the step size $h$, the harmonic average finite difference equation at a grid point $x_{i}$ is the following
\eql{hav}
 \frac{1}{h^2}  \left( \beta _{ i+ \frac{1}{2} } { ( U_{i+1} -U_i)} -
\beta _{i- \frac{1}{2} } { ( U_i -U_{i-1} )} \right)  U_i - \sigma_i  \, U_i= f_i,
\en
 where  $\sigma_i = \sigma(x_i), \; \;f_i = f(x_i)$, and $\beta_{i+\half}$ is the following harmonic average,
\eqm
  \beta _{ i+ \frac{1}{2} } = \left (  \frac{1}{h}  \int_{x_i}^{x_{i+1}} \beta^{-1}(x) dx   \right)^{-1}.
 \enm
 When $\beta(x)\in C[0,1]$,  the above expression can be replaced with $\beta(x_{ i+ 1/2})$. 
 Numerically,   second order accuracy in the pointwise $L^{\infty}$ norm has been observed even if $\beta$ has a finite jump discontinuity but the solution and the flux are continuous, that is, $w=0$ and $v=0$ assuming that the integration is accurate enough, see for example \cite{rutka-phd-thesis,li:thesis}. Note that the finite difference scheme is inconsistent at the two grid points adjacent to the discontinuity as we can see later. One misconception is that the local truncation errors would cancel out to $O(h)$,  which would lead to second order convergence. The harmonic average finite difference method does not converge when the solution or flux has a finite jump, which is often referred as a non-homogeneous jump condition.

In this paper, we not only show the convergence of the harmonic average finite difference method for 1D interface problems with finite number of discontinuities in the coefficient with homogeneous jump conditions ($w=0$ and $v=0$), but also proposed a new improved second order harmonic average finite difference method for the general 1D interface problem \eqn{eq1dB} with non-homogeneous jump conditions.

The paper is organized as follows. In the next section, we explain the classical and the new improved harmonic average finite difference method. In Section~\ref{sec:conv}, we show the convergence analysis of the classical and the new improved harmonic average method.  In Section~\ref{sec:2D}, we present the improved harmonic average method for two-dimensional problems when the interface is parallel to one of axis. Numerical experimental results are shown in Section~\ref{sec:examples}. We conclude in the last section. 

\section{The improved second order harmonic average finite difference method} \label{sec:HAV}

\begin{figure}[!htbp]
  \centerline{\includegraphics[width=0.45\textwidth]{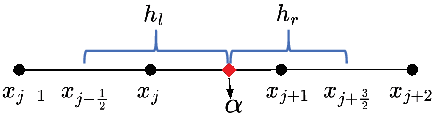}}
 \caption{A diagram of a one-dimensional grid and the interface $\alpha$ between $x_j$
and $x_{j+1}$.} \label{domain}
\end{figure}

The improved second order harmonic average finite difference method with non-homogeneous jump conditions and $\sigma(x)=0$\footnote{If $\sigma(x)\neq 0$ but in $C^0$, we can treat it as a source term. So, the method and analysis are still valid.} is outlined below.
\bi
\item {\bf Step 1:} Generate a Cartesian grid, say
\eqmno
  x_{i} = i h, \;\;i = 1,2, \cdots, N,
\enmno
where $ h = 1/N.$ The point $\alf $ will typically fall between two
grid points, say $ x_{j} \leq \alf < x_{j+1}$ (see Fig. \ref{domain}). The grid points $ x_{j}$ and
$x_{j+1}$ are called {\em irregular grid points}
if a standard three-point central
finite difference stencil is going to be used at grid points away
from the interface $\alf$. The other grid points
are called regular grid points.

\item {\bf Step 2:} Determine the finite difference equation at regular grid points.
At a grid point $x_{i}$, $ i \neq j$, $j+1$,
the standard three-point central finite difference approximation
\eql{FD_regular}
 \frac{1}{h^2}  \left( \beta _{ i+ \frac{1}{2} } { ( U_{i+1} -U_i)} -
\beta _{i- \frac{1}{2} } { ( U_i -U_{i-1} )} \right)  = f_i,
\en
is used, where $\beta_{i+\half} = \beta(x_{{i+\half}} )$ and so on, $
\sigma_i = \sigma(x_i), \; \;f_i = f(x_i)$.

\item {\bf Step 3:} Determine the finite difference scheme at irregular regular grid
points $x_j$ and $x_{j +1}$. The finite difference equations
are:   
\eqml{hav_FD}
& \dsp \frac{1}{  h_l}  \left(   \frac{ \bar{\beta}_{j+\frac{1}{2}}  ( U_{j+1} -U_j)}{h }   -
\frac{\ \beta _{j- \frac{1}{2} } ( U_j -U_{j-1} )}{h}    \right)   = f_j + C_j, \\  \eqsp \eqsp
& \!\!\!  \!\!\!  \!\!\! \dsp  \frac{1}{ h_r}  \left(  \frac{ \beta_{ j+ \frac{3}{2} } ( U_{j+2} -U_{j+1}  ) }{h} -
 \frac{ \bar{\beta}_{j+\frac{1}{2}} ( U_{j+1} -U_j)}{h }
 \right)    =  f_{j+1} + C_{j+1},
\enml
where $h_l= \alf - x_{j-\half}>0$, $h_r =  x_{j+\frac{3}{2}}-\alf >0$, and
\eqm
   \bar{\beta}_{j+\frac{1}{2}}  & = & \dsp \left (  \frac{ x_{j+1}-\alf} { \beta^+ h} +   \frac{ \alf-x_j} { \beta^- h}  \right)^{-1}, \label{eqbeta} \\ \eqsp
  C_j &= & \frac{ \bar{\beta}_{j+\frac{1}{2}} }{h_l\,h} \left ( \frac{\null}{\null} [u] +  {\frac{[\beta u_x]}{\beta^+}} (x_{j+1}-\alf ) \right ), \label{eqcj} \\ \eqsp
 C_{j+1} &= &- \frac{ \bar{\beta}_{j+\frac{1}{2}} }{h_r \, h  } \left ( \frac{\null}{\null} [u] + {\frac{[\beta u_x]}{\beta^-}} (x_{j}-\alf ) \right ) .
 \enm
\item {\bf Step 4:} Solve the system of equations \eqn{FD_regular}-\eqn{hav_FD}, denoted as $A_h \,\mathbf{U} =\mathbf{F}$,
whose coefficient matrix is tridiagonal, to
get an approximate solution of~$u(x)$ at all grid points.
\ei
 We call the finite difference scheme above as the {\em improved harmonic average method}. 

\vthin

The  improved harmonic average method is similar to the traditional method for non-interface problems and does not need the derivative of $\beta(x)$ which is different from other methods.
We see also that the coefficient matrix of the finite difference equations is diagonally dominant, and the off-diagonals have the opposite sign as that of the main diagonals and it is irreducible. If one of boundary conditions is prescribed (Dirichlet), then the coefficient matrix $A$ is an M-matrix and invertible. The key to the improved harmonic average method is the added correction terms $C_j$ and $C_{j+1}$ corresponding to the jumps in the solution and the flux. 

\begin{remark} In fact, the harmonic average in the interval $(x_j, x_{j+1})$
\eqml{new_beta_j}
     \beta _{j+ 1/2 }  &=&\dsp  \left (  \frac{1}{h}  \int_{x_j}^{x_{j+1}} \beta^{-1}(x) dx   \right)^{-1} = \left (  \frac{1}{h}  \int_{x_j}^{\alpha} \beta^{-1}(x) dx  + \frac{1}{h}  \int_{\alpha}^{x_{j+1}} \beta^{-1}(x) dx  \right)^{-1} \\ \eqsp 
    &=& \dsp  \left (  \frac{1}{h}  \int_{x_j}^{\alpha} \frac{1}{\beta^- + O(h)} dx  + \frac{1}{h}  \int_{\alpha}^{x_{j+1}} \frac{1}{\beta^+ + O(h)} dx  \right)^{-1} \\ \eqsp
    & = & \dsp \left (  \frac{1}{h}  \int_{x_j}^{\alpha} \frac{1}{\beta^-} dx  + \frac{1}{h}  \int_{\alpha}^{x_{j+1}} \frac{1}{\beta^+} dx  + O(h)\right)^{-1}     \\ \eqsp
    & =&  \dsp \bar{\beta}_{j+\frac{1}{2}} + O(h).
\enml
When the coefficient $\beta(x)$ is a piecewise constant, the harmonic average $\beta _{j+ \frac{1}{2}} = \bar{\beta}_{j+\frac{1}{2}}$.
\end{remark}

\begin{remark}
In general, $A$ is not symmetric anymore except that $h_l=h_r$, that is, the interface is in the middle of two grid points.
However, we can scale the finite difference equations at the two grid points to make the coefficient matrix symmetric by multiplying $h_l/h$ and $h_r/h$ to the first and second equations in \eqn{hav_FD}, respectively.
\end{remark}

\begin{remark}
  For one dimensional interface problems with a {\em variable} discontinuous coefficient, the immersed interface method proposed in \cite{li:book,li:thesis} can work but needs the first order derivative of the coefficient from the left and right. The improved harmonic average method is derivative free, and the coefficient matrix of the finite difference equations is an M-matrix, which is a big advantage. The immersed finite element method (IFEM), see for example, \cite{li:fem,li:tao1,Lin_Lin_Zhang_SINUM15,Ji_Li_wang_chen} is a derivative free method, but the convergence proof is only valid in the $H^1$ or $L^2$ norms, not in the $L^{\infty}$ norm.
\end{remark}

\begin{remark}
  If $\sigma(x) $ in \eqn{eq1d} is non-zero but continuous, then we can treat $\sigma(x) u(x)$ term as a source term, which is $\sigma(x_i) U_i$ in the discrete case added to the finite difference equations. The finite difference scheme can be further improved with an additional correction term if the source term $f(x)$ has a finite jump at $\alf$, see for example \cite{rjl-li:sinum,li:book}.
\end{remark}

\subsection{The derivation of the  improved harmonic average method   at
an irregular grid point}

The local truncation errors of the finite difference scheme  are defined as,
\eqm
T_i = \frac{1}{h^2}  \left( \beta _{ i+ \frac{1}{2} } { ( u(x_{i+1} )-u(x_i))} -
\beta _{i- \frac{1}{2} } { ( u(x_i)  -u(x_{i-1} ) )} \right)   - f_i,
\enm
for $i=1,2,\cdots, x_{n-1}$,  $i\neq j, \, j+1$.  At $x_j$ and $x_{j+1}$, the local truncation errors are,
\eqm
T_j = \frac{1}{  h_l}  \left(  \frac{ \bar{\beta}_{j+\frac{1}{2}}  ( u(x_{j+1}) -u(x_j) )} {h }   -
\frac{\ \beta _{j- \frac{1}{2} } ( u(x_j) -u(x_{j-1})  )}{h}    \right)    - f_j -C_j,
\enm
\eqml{Tj_j+1}
 T_{j+1} &=& \dsp \frac{1}{ h_r}  \left(  \frac{ \beta_{ j+ \frac{1}{2} } ( u(x_{j+2}) -u(x_{j+1} ) ) }{h} -
 \frac{   \bar{\beta}_{j+\frac{1}{2}}  ( u(x_{j+1} ) -u(x_j) ) }{h  }  \right)   \\ \eqsp
 && \quad \qquad  - f_{j+1} \textcolor{blue}{-} C_{j+1}. 
\enml

At regular grid points, we know that $|T_i |\le C h^2$. We also know that $(u(x_{j+1}) -u(x_j) )/h= u_x(x_{j+\frac{1}{2}})+ O(h^2)$.
At irregular grid points $x_j$ and $x_{j+1}$, we use the IIM's idea and an underdetermined coefficients method to show the finite difference weights and the correction terms. We use the irregular grid point $x_j$ to explain the process.

 Using the jump conditions and Taylor expansions at $x=\alf$ we derive
\eqmno
  \beta^- u_x^- &\approx& \bar{\beta}_{j+\frac{1}{2}} \frac{ u(x_{j+1}) -u(x_j) }{h} - \bar{C}_j  \\ \eqsp
  &=& \frac{\bar{\beta}_{j+\frac{1}{2}}}{h} \left ( u^+ + u_x^+ ( x_{j+1}-\alf)   - u^- - u_x^-(x_j-\alf) + O(h^2) \right ) - \bar{C}_j \\ \eqsp
   &=& \frac{\bar{\beta}_{j+\frac{1}{2}}}{h} \left ( u^- + [u] +  \frac{ \beta^- u_x^- + [\beta u_x] }{ \beta^+} ( x_{j+1}-\alf) - u^- - u_x^-(x_j-\alf) + O(h^2) \right ) - \bar{C}_j \\ \eqsp
   &=& \frac{\bar{\beta}_{j+\frac{1}{2}}}{h} \left ( \Big( \frac{ \beta^-   }{ \beta^+} ( x_{j+1}-\alf) + (\alf - x_j) \Big)u_x^- +  [u] +  \frac{ [\beta u_x] }{ \beta^+} ( x_{j+1}-\alf)  \right ) - \bar{C}_j + O(h),
\enmno
where again $u_x^{\pm}$ is the limit of $u'(x)$ from the right and the left of $\alpha$. By matching the  corresponding terms on both sides, we  should set
\eqmno
   &  \dsp \frac{\bar{\beta}_{j+\frac{1}{2}}}{h}  \left (  \frac{ \beta^-   }{ \beta^+} ( x_{j+1}-\alf) + (\alf - x_j)  \right ) = \beta^-, \\ \eqsp
&\dsp  \bar{C}_j = \frac{\bar{\beta}_{j+\frac{1}{2}}}{h}  \left (  [u] +  \frac{ [\beta u_x] }{ \beta^+} ( x_{j+1}-\alf)  \right ).
 \enmno
Solving the above two linear equations, we obtain the formula (\ref{eqbeta}) for the underdetermined coefficient $\bar{\beta}_{j+\frac{1}{2}}$ and the correction term $\bar{C}_j$ at  $x_j$,
\eqmno
  \bar{C}_j = C_jh_l.
\enmno

  After we obtain the approximate to $\beta^- u_x^-$, we have,
\eqmno
 \left.  \frac{\null}{\null} \left (\beta(x) u_x \right )_x \right |_{x=x_j} &\approx& \frac{(\beta(x) u_x)|_{\alpha^-} - (\beta(x) u_x)|_{x_{j-\frac{1}{2}}}}{\alpha - x_{j-\frac{1}{2}}} \\
  && \approx \frac{\bar{\beta}_{j+\frac{1}{2}} ( U_{j+1} -U_j ) - \bar{C}_j h - \beta_{j- \frac{1}{2} } ( U_j -U_{j-1} ) }{h h_l}  \\
  & &= \frac{1}{  h_l}  \left(   \frac{ \bar{\beta}_{j+\frac{1}{2}}  ( U_{j+1} -U_j)}{h }   -
\frac{\ \beta _{j- \frac{1}{2} } ( U_j -U_{j-1} )}{h}    \right)  - C_j.
\enmno
  Using the above expression to approximate the ODE (\ref{eq1d}), we obtain the finite difference equation (\ref{hav_FD}) at the irregular grid point $x_j$. The derivation for the grid point $x_{j+1}$  is similar, so is omitted.

 In the next section,  we derive the estimates of the local truncation errors $T_j$ and $T_{j+1}$.

\section{Theoretical analysis of the classical and improved harmonic finite difference scheme} \label{sec:conv}

\begin{lemma}
 Under the assumptions of the interface problem \eqn{eq1d} and the algorithm settings, and assuming that the solution to \eqn{eq1d} is piecewise $C^3(0,1)$, then the  following hold:
 \eqm
  T_j  &=&   \frac{   \bar{\beta}_{j+\frac{1}{2}} } {h h_l} \frac{ u^+_{xx} (x_{j+1} - \alf)^2 - u^-_{xx} (x_j - \alf)^2 }{2} + O(h) , \label{trunj}  \\ \eqsp
   T_{j+1} &=& - \frac{  \bar{\beta}_{j+\frac{1}{2}}} { h \,h_r} \frac{ u^+_{xx} (x_{j+1} - \alf)^2 - u^-_{xx} (x_j - \alf)^2 }{2} + O(h). \label{trunjr}
 \enm
\end{lemma}
Note that, we can see clearly that $T_j\sim O(1)$ and $T_{j+1}\sim O(1)$ in general, which indicates the algorithm is inconsistent. It is easy to spot that,
\eqm
   T_j + T_{j+1}\sim O(1),  \quad \mbox{but} \quad T_{j+1} = - \frac{h_r}{h_l} \, T_j  \label{Tj_Tj+1} + O(h).
\enm
If $h_r = h_l$, then $T_j = - T_{j+1} + O(h)$,  we would have the error cancellation, which does not happen in general. 
  The $O(1)$ terms contains $(u^+_{xx} (x_{j+1} - \alf)^2 - u^-_{xx} (x_j - \alf)^2)/(h_l h  )$ and so on, which fluctuates depending on the relative position of $\alf$ between $x_j$ and $x_{j+1}$. So, there are some cancellations in the local truncation errors but not enough to conclude second order convergency.  We need additional tools to show that $A_h^{-1} (T_j{\bf e}_j  + T_{j+1}{\bf e}_{j+1}  ) \sim O(h^2)$, which leads to the second order convergence.
To prove the convergence of the finite difference scheme (improved harmonic average method), we first prove the piecewise constant coefficient $\beta$ that has a finite jump at $\alf$. For this purpose,
we introduce the following second-kind Green's function,  see the similar one in \cite{IIM_Deriv18}.

\vthin

  Let $G_{\alpha}(x)$ be  the solution of the following interface problem,
\eqml{1dgreen}
  && (\beta G_{\alpha}'(x))'  = \delta'(x-\alf) , \qquad 0 < x < 1, \\ \eqsp
  && G_{\alpha}(0)  = 0, \quad G_{\alpha}(1)  = 0.
\enml
Note that there is a  jump at $x= \alf$ for the solution but not the flux, that is, 
\eqm
  \left[ \frac{\null}{\null} G_{\alpha}(x )\, \right]_{\alpha} = \frac{2} {\beta^-+\beta^+}, \qquad \left[ \frac{\null}{\null} \beta \, G_{\alpha}'(x)\, \right]_{\alpha} = 0,
\enm
see \cite{REU10}   It is easy to check that the interface problem has  the following analytic solution,
\begin{eqnarray} \label{green}
G_{\alpha} (x )=
\begin{cases}
\dsp  \left(   \frac{1 }{\beta^-(1-\alpha ) + \beta^+ \alpha} \right ) \frac{-2\beta ^+}{\beta^- + \beta^+}  \, x ,& 0< x< \alpha, \\ \eqsp
 \dsp   \left(    \frac{1 }{\beta^-(1-\alpha ) + \beta^+ \alpha}  \right ) \frac{2\beta ^-}{\beta^- + \beta^+} \, (1-x),& \alpha < x < 1,
  \end{cases}
\end{eqnarray}
where  $\beta^-(1-\alpha ) + \beta^+ \alpha>0$.

If we apply the improved harmonic average method  to the 1D interface problem \eqn{1dgreen}, we will get the exact solution (without the presence of the round-off errors) since the second order derivatives are zero excluding the interface $\alpha$, which implies that all the local truncation errors are zero. Let us write   the finite difference equations as a matrix-vector form $A_h {\bf G} = {\bf \tilde{F}}$.  Then, the $j$-th component of $ {\bf \tilde{F}}$ is 
$$\tilde{F}_j=   \frac{ \bar{\beta}_{j+\frac{1}{2}} }{h_l\,h} \frac{2}{\beta^- + \beta^+},$$ 
from the correction term and the exact solution. Thus, if we take $W = T_j h_l h (\beta^- + \beta^+) /  (2\bar{\beta}_{j+\frac{1}{2}}) $, then we have $A_h (G_{\alpha} (x_j )W) =T_j$.
We summarize the result in the following lemma.

\begin{lemma}
 With the settings in this section, we have
\begin{equation}\label{aa}
 \left  (A_h^{-1}  \left (T_j \,{\bf e}_j  + T_{j+1} \,{\bf e}_{j+1}  \frac{\null}{\null} \right ) \right )_i =
G_\alpha(x_i) W,  
\end{equation}
 where ${\bf e}_j$ is the $j$-th unit vector whose entries are all zero  except for the $j$-th entry that is number one, and $W =  T_j h_l h (\beta^- + \beta^+) /  (2\bar{\beta}_{j+\frac{1}{2}})$. It is obvious that $|W| \le C h^2$.
 \end{lemma}

\noindent\emph{Proof:} First we note that $G_{\alpha} (x )$ is a piecewise linear function and we would obtain the exact solution if we apply the algorithm to the boundary value problem,
\eqml{eq1dB}
& \dsp (\beta(x) u_x)_x  = 0, \qquad \{0 < x < \alf\} \cup \{  \alf < x <1\}, \\ \eqsp
& \dsp [u]_{\alf} = \frac{2}{\beta^- + \beta^+} W, \qquad [\beta u_x] _{\alf} = 0, \qquad 0 < \alf < 1, \\ \eqsp
& \dsp u(0) = 0, \quad u(1) = 0,
\enml
which is equivalent to \eqn{1dgreen}. 
At the grid point $x_j$, the right hand side is simply the correction term,
\eqm
  {C}_j = \frac{   \bar{\beta}_{j+\frac{1}{2}} }{h_l \, h  } \, \frac{2}{\beta^- + \beta^+} W  = T_j,
\enm
and at the grid point $x_{j+1}$, the right hand side is
\eqm
  {C}_{j+1} = -\frac{\bar{\beta}_{j+\frac{1}{2}}}{h_r \, h } \, \frac{2}{\beta^- + \beta^+} W =  -\frac{ h_l}{h_r } \, T_j  = T_{j+1},
\enm
from \eqn{Tj_j+1}.
The rest of the right hand sides are all zero. Since the algorithm would return the true solution if it is apply it to \eqn{eq1dB},  we conclude (\ref{aa}) and
\eqml{Ejj+1}
  |\left ( A_h^{-1} \left ( T_j \mathbf{e}_j+  T_{j+1} \mathbf{e}_{j+1} \right ) \right )_i| &= &  \dsp |G_\alpha(x_i)| \cdot |W| \le C h^2.
\enml

Now we state the main convergence theorem on the convergence of the improved harmonic average method that includes the original harmonic average method as a special case when $[u]=0$ and $[\beta u_x]=0$. 

\begin{theorem} \label{theo1d}
With the assumptions and setting is this section,  the finite difference solution  computed using the improved harmonic average method applied to \eqn{eq1d} or \eqn{eq1dB}  has second order convergence in the infinity norm, that is,
\eqm
  \|\mathbf{E} \|_{\infty}\le C h^2,
\enm
assuming that the true solution of \eqn{eq1d} is piecewise $C^4$ excluding the interface $\alf$, that is, 
$u(x) \in C^4(0,\alf)  \cup C^4(\alf,1)$. 
\end{theorem}

\noindent\emph{Proof:} We have
\eqm
  A_h \mathbf{U} = \mathbf{F}, \qquad A_h \mathbf{u} = \mathbf{F} + \mathbf{T}, \qquad \mbox{thus} \quad A_h \mathbf{E} = \mathbf{T},
\enm
where $\mathbf{T}$ is the local truncation error vector with $|T_i|\le C h^2$ at regular grid points, and $T_j $, and $T_{j+1}$ are given in \eqn{trunj} and \eqn{trunjr}.  Thus, we get
\eqmno
    E_i &=& \sum_{k = 1}^{N - 1}\left ( A_h^{-1}  T_k  \mathbf{e}_k \right )_i \\
    & &= \sum_{k \neq j,j+1} \left ( A_h^{-1}  T_k  \mathbf{e}_k \right )_i + \left ( A_h^{-1} \left ( T_j \mathbf{e}_j+  T_{j+1} \mathbf{e}_{j+1} \right ) \right )_i \\
    & &=  O(h^2) + \left ( A_h^{-1} \left ( T_j \mathbf{e}_j+  T_{j+1} \mathbf{e}_{j+1} \right ) \right )_i.
\enmno
From \eqn{Ejj+1}, we have $|E_i| \le Ch^2$, which completes the proof of the theorem.

\begin{remark} $\null$

\be

  \item  In order to have second order accuracy, we still need that condition that $u(x)$ is piecewise $C^4$ so that the local truncations errors are $O(h^2)$ at regular grid points although we just need $C^3$ condition near the interface.
  
  \item As we can see from the proof, the error constant, while it is $O(1)$, depends on $h_l$ and $h_r$, which are mesh dependent quantities. Also
the leading errors at irregular grid points will depend on $(x_{j}-\alf)^3$ and $(x_{j+1}-\alf)^3$, so the convergence order may fluctuate as confirmed from the numerical experiments.

  \item While the proof is for piecewise constant coefficients, the convergence analysis is still valid for variable coefficients since we can use piecewise Taylor expansion at the interface for the coefficient. After the expansion, the $O(h)$ term can contribute at most $O(h^2)$ to the global error.
  \item For the equation $(\beta(x) u_x)_x - \sigma(x) u = f(x)$, we can treat $\sigma(x) u(x)$ as a source term and the developed method can still be applied if $\sigma(x)$ is continuous at $\alpha$.
\ee
\end{remark}

\section{The improved harmonic method for two dimensional problems with an interface that is parallel to one of axis}
\label{sec:2D}

The improved harmonic method can be applied directly to two- or three-dimensional interface problem with a line ( plane in 3D) interface that is parallel to one of axis. Consider for example,
\eqml{eq2dB}
& \dsp (\beta(x,y) u_x)_x + (\beta(x,y) u_y)_y - \sigma(x,y) u = f(x,y), \quad
\begin{array}{l}
     \{a < x < \alf\} \cup \{  \alf < x <b\},   \\ \eqsp
    \{  c < y <d\} ,
  \end{array} \\ \eqsp
& \dsp [u]_{\alf} = w, \qquad [\beta u_x] _{\alf} = v.
\enml
We assume that $\sigma(x,y) $ and $f(x,y)$ are continuous in the domain, and a Dirichlet boundary condition along the boundary.  We assume that the coefficient $\beta(x,y)$ has a finite jump discontinuity in the $x$ direction at $x=\alf$.
After we set up a uniform grid
\eqm
 x_i = a + i h_x, \quad i=0,1,\cdots, m;  \quad y_j = c + i h_y, \quad j=0,1,\cdots, n-1,
\enm
with $h_x = (b-a)/m$; $h_y= (d-c)/n$.  Let us assume that $x_k \le \alf <x_{k+1}$, which was $x_j$ before for the 1D problem. The finite difference equation at a grid point $(x_i,y_j)$ is
\eqml{hav_FD2-2d}
&& \dsp \frac{1}{h_x^2}  \left( \beta _{ i+ \frac{1}{2},j } { ( U_{i+1,j} -U_{ij})} -
\beta _{i- \frac{1}{2},j } { ( U_{ij} -U_{i-1,j} )} \right)   \\ \eqsp
&& \dsp \qquad \null + \frac{1}{h_y^2}  \left( \beta _{ i,j + \frac{1}{2}} { ( U_{i,j+1} -U_{ij})} -
\beta _{i,j - \frac{1}{2}} { ( U_{ij} -U_{i,j-1} )} \right)   -\sigma_{ij} U_{ij} = f_{ij}  ,
\enml
for $i=1,2,\cdots, m$, $i\neq k, k+1$, $j=1,2,\cdots, n-1$. At the irregular grid points $(x_k,y_j)$, the finite difference equations are
\eqml{hav_FD2-2dj}
&& \dsp \frac{1}{  h_l}  \left(   \frac{ \bar{\beta}_{k+\frac{1}{2},j}  ( U_{k+1,j} -U_{kj})}{h_x }   -
\frac{\ \beta _{k- \frac{1}{2},j } ( U_{kj} -U_{k-1,j} )}{h_x}    \right)   \\ \eqsp
&& \dsp \qquad \null  + \frac{1}{h_y^2}  \left( \beta _{ k,j + \frac{1}{2}} { ( U_{k,j+1} -U_{kj})} -
\beta _{k,j - \frac{1}{2}} { ( U_{kj} -U_{k,j-1} )} \right)   -\sigma_{kj} U_{kj}  = f_{kj} + C_{kj}, \\  \eqsp \eqsp
&& \dsp C_{kj} = \frac{\bar{\beta}_{k+\frac{1}{2},j}}{h_l  h_x}  \left (  [u]_{\alf,y_j}  +  \frac{ [\beta u_x]_{\alf,y_j} }{ \beta^+} ( x_{k+1}-\alf)  \right ),
\enml
 $j=1,2,\cdots, n$,  $h_l= \alf - x_{k-\half}>0$, and
 \eqm
  \bar{\beta}_{k+\frac{1}{2},j}  & = & \dsp \left (  \frac{ x_{k+1}-\alf} { \beta^+_{\alf,y_j}\,  h_x} +
  \frac{ \alf-x_k} { \beta^-_{\alf,y_j} \,h_x}  \right)^{-1},
 \enm
 have the similar meanings as defined before.

Similarly, at the irregular grid points $x_k$, the finite difference equations are
\eqml{hav_FD2-2djp}
&& \dsp \frac{1}{ h_r}  \left(  \frac{ \beta_{ k+ \frac{3}{2},j } ( U_{k+2,j} -U_{k+1,j}  ) }{h_x} -
 \frac{ \bar{\beta}_{k+\frac{1}{2},j} ( U_{k+2,j} -U_{k+1,j} )}{h_x }
 \right)    \\ \eqsp
&& \dsp \qquad \null  + \frac{1}{h_y^2}  \left( \beta _{ k+1,j + \frac{1}{2}} { ( U_{k+1,j+1} -U_{k+1,j})} -
\beta _{k+1,j - \frac{1}{2}} { ( U_{k+1,j} -U_{k+1,j-1} )} \right)   \\ \eqsp
&& \null \qquad \qquad \qquad \dsp -\sigma_{k+1,j} U_{k+1,j}  = f_{k+1,j} + C_{k+1,j}, \\  \eqsp \eqsp
&& \dsp C_{k+1,j} =- \frac{ \bar{\beta}_{k+\frac{1}{2} ,j} }{h_r \, h _x } \left ( \frac{\null}{\null} [u]_{\alf,y_j} +  \frac{[\beta u_x]_{\alf,y_j}}{\beta^-} (x_{k}-\alf ) \right ) ,
\enml
 $j=1,2,\cdots, n$, and $\bar{\beta}_{k+\frac{1}{2},j}$ and $h_r =  x_{k+\frac{3}{2}}-\alf >0$ have the same meanings as defined before.

\section{Numerical examples} \label{sec:examples}

We validate the method and analysis through some nontrivial examples for which we have analytic solutions.

\begin{example}
 Let the coefficients  and the source term be the following,
 \eqm
  \beta(x) & = & \left \{ \begin{array}{ll}
1+ x^2     &   \mbox{if    $0<x< \alf $,} \\ \eqsp
  \log(2 + x)    &    \mbox{ if $ \alf < x < 1; $}
  \end{array} \right. \\ \eqsp
  f(x) & = & \left \{ \begin{array}{ll}
2 x k_1 \cos  (k_1 x)   - (x^2 + 1) k_1^2 \sin(k_1 x)     &   \mbox{if    $0<x< \alf $,} \\ \eqsp
\dsp   -\frac{k_2 \sin  (k_2 x)  }{ 2 + x} - \log(x+2) k_2^2  \cos (k_2 x)     &    \mbox{ if $ \alf < x < 1. $}
  \end{array} \right.
  \enm
The solution to the boundary value problem is
\eqm
 u(x)  &=& \left \{ \begin{array}{ll}
 \sin(k_1 x)     &   \mbox{if    $0<x< \alf $,} \\ \eqsp
  \cos (k_2 x)    &    \mbox{ if $ \alf < x < 1. $}
  \end{array} \right.
\enm
\end{example}

In Table~\ref{tab1}, we show a grid refinement analysis of the method with  $\alf=1/3$, and different $k_1$ and $k_2$ for the example.
The first column is the mesh size $n$. The second to fifth columns are results of the solution errors, the local truncation errors at the left and the right irregular grid points when $k_1 = 5, k_2 =3$. We clearly say that the local truncation errors at the two grid points are $O(1)$. The approximated convergence order computed by
\eqm
  \mbox{order} = \frac{\log ( \|E \|_{\infty, n}/\|E \|_{\infty, 2n})} {\log 2} .
\enm
The convergence order jumps up and down, but the average is $2.1093$ since the error constant depends on the relative position of the interface and the underline grid. In fact, the averages of the orders of the errors using two consecutive grids ($n$ and then $2n$) is close the number $2$.  The last two columns in Table~\ref{tab1} are the results of the method with  $\alf=1/3$ with $k_1=30, k_2 = 25$. We see larger errors but the same behavior of the convergence orders. The local truncation errors at the two grid points are propositional to $k_1^2 \sim 900/2$ from the second order derivatives. The actual local truncation errors at the two irregular grid points vary between $5$ and $112$.

The numerical example in Table~\ref{tab1} also confirmed our theoretical analysis. While the harmonic average method works when the coefficient has finite jumps, it is inconsistent because the local truncations at two grid points neighboring the interface, say $T_j$ and $T_{j+1}$ are of $O(1)$. Secondly, let $\hat{{\bf T}}_j$ be the vector whose components are all zero except the $j$-th component which  $T_j$, the local truncation error at $x_j$, the irregular grid point from the lest, then $A_h^{-1} \hat{{\bf T}}_j \sim O(h)$, so is  $A_h^{-1} \hat{{\bf T}}_{j+1} \sim O(h)$. Nevertheless,  $A_h^{-1} \left( \hat{{\bf T}}_j + \hat{{\bf T}}_{j+1} \right )\sim O(h^2)$, which tells us how the errors cancellations take place.

\begin{table}[hptb]
\centering \caption{Grid refinement analysis of the method with  $\alf=1/3$, and different $k_1$ and $k_2$.
} \label{tab1}

\vthin

\begin{tabular}{||c||c|c||c| c||c| c||}
\hline
     $N$ &   $\|E \|_{\infty}$ & order  & $T_j $  &   $T_{j+1} $    & $\|E \|_{\infty}$  & order \\  \hline
     & \multicolumn{4} {|c||}{$k_1 = 5, k_2 =3$}  &  \multicolumn{2} {|c|}{$k_1 = 30, k_2 =25$} \\ \hline

32 &   2.7133e-03  &  &  4.5144   &   -6.2990   & 	  8.5415e-01  & \\   \hline
64 &   4.0582e-04   & 2.7411 & 3.8603e-01  & -2.6182e-01  & 1.4703e-01 & 2.5384  \\   \hline
128 &  1.8006e-04   &   1.1724 & 4.5277     & -6.3423    &   5.3686e-02 & 1.4535\\   \hline
256 &  2.5043e-05   &2.8460  & 3.4366e-01 &  -2.4413e-01 &  9.1750e-03 & 2.5488 \\   \hline
512 &  1.1418e-05   & 1.1331 & 4.5354   &   -6.3510    &   3.3278e-03  & 1.4631 \\   \hline
1024  &  1.5603e-06  & 2.8714 &  3.3574e-01  & -2.3962e-01  &  5.7223e-04 & 2.5399\\   \hline
2048 &   7.1622e-07   & 1.1233 & 4.5376    &   -6.3531   &    2.0742e-04  & 1.4640\\   \hline
4096 & 9.7443e-08  &2.8778 &  3.3394e-01  & -2.3848e-01  &  3.5742e-05 & 2.5369  \\ \hline
 & average & 2.1093 &O(1) & O(1)&average & 2.0778  \\ \hline
   \end{tabular}
\end{table}

\begin{figure}[hpt]
\begin{minipage}[t]{3.0in} (a)\\

\includegraphics[width=0.85\textwidth]{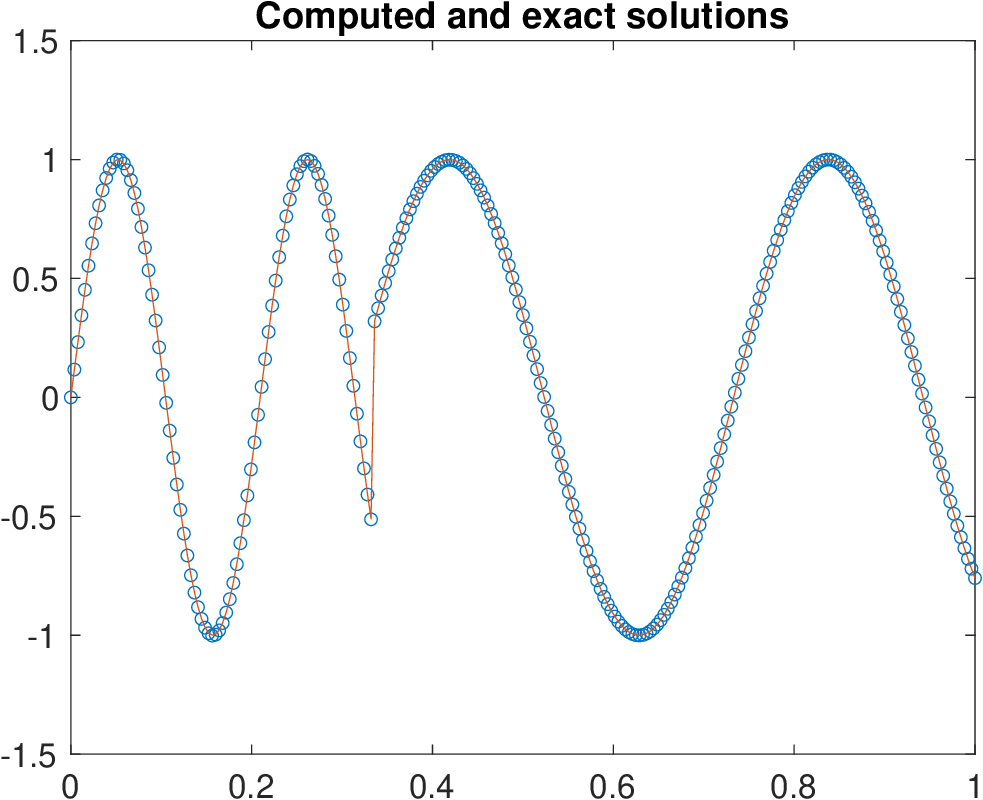}
 \end{minipage}
\begin{minipage}[t]{3.0in} (b)\\

\includegraphics[width=0.88\textwidth]{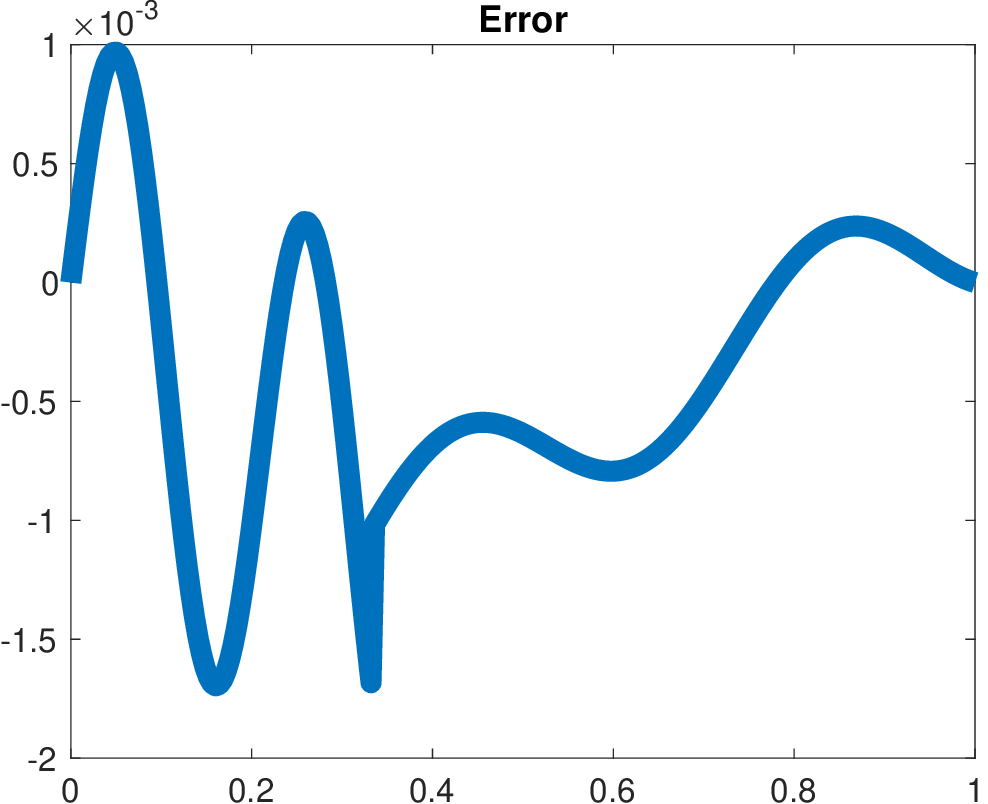}
 \end{minipage}

\hfill

\caption{(a), computed (little 'o's) and exact solution (solid line) when $k_1=30, k_2=5$ and $n=256$.  (b), the error plot. }\label{fig:Ex1}
\end{figure}

In Figure~\ref{fig:Ex1}~(a), we show the computed solution (dotted line) and the exact solution (solid line) with $\alpha = 1/3$, $k_1=30, k_2=5$ and $n=256$, which are  identical by naked eyes although the solution and the flux are discontinuous. In Figure~\ref{fig:Ex1}~(b), we plot the errors of the computed solution. Note that the errors near the interface ($1/3$) have roughly the same magnitude as those at regular grid points  (around $0.18$).

\begin{example} $\null$
\end{example}
Next. we consider the  effect of the jump ratio $\beta^-/\beta^+$ on the algorithm. To do so, we change the setting as below,
 \eqm
  \beta(x) & = & \left \{ \begin{array}{ll}
\beta^-   &   \mbox{if    $0<x< \alf $,} \\ \eqsp
  \beta^+    &    \mbox{ if $ \alf < x < 1; $}
  \end{array} \right. \\ \eqsp
  f(x) & = & \left \{ \begin{array}{ll}
   -   \beta^-   k_1^2 \sin(k_1 x)     &   \mbox{if    $0<x< \alf $,} \\ \eqsp
\dsp   - \beta^+   k_2^2  \cos (k_2 x)     &    \mbox{ if $ \alf < x < 1, $}
  \end{array} \right.
  \enm
while the solution is unchanged. In Table~\ref{tab2}, we list the grid refinement results for different jump ratios. We still observe average second order convergence for all different jump ratios (small $1.5/3$ and  large $2000/1.5$ or $1.5/2000$) although the errors can increase slightly at some mesh sizes followed by significant decrease in the next. So, the method seems to be insensitive to the jump ratios.

\begin{table}[hptb]
\centering \caption{Grid refinement analysis of the method with  $\alf=1/3$,  $k_1=5$ and $k_2=3$,  and different jump ratios.
} \label{tab2}

\vthin

\begin{tabular}{||c||c|c||c| c||c| c||}
\hline
     $N$ &   $\|E \|_{\infty}$ & order  & $\|E \|_{\infty}$ & order  & $\|E \|_{\infty}$ & order   \\  \hline
 & \multicolumn{2} {|c||}{$\beta^-=1.5, \beta^+ = 3$}  &  \multicolumn{2} {|c||}{$\beta^-=2000, \beta^+ = 1.5$}  & \multicolumn{2} {|c|}{$\beta^-=1.5, \beta^+ = 2000$}  \\ \hline
  32 &4.3434e-03      &        &  2.7967e-03       &        & 6.3934e-03       &       \\  \hline
   64 & 2.3005e-04   & 4.2388   &  4.9463e-04  & 2.4993  &  2.3391e-04   & 4.7726   \\  \hline
  128 &  2.8438e-04 & -0.3058  & 1.6432e-04  & 1.5898  &  4.1245e-04  & -0.8182 \\  \hline
256 &   1.4146e-05   & 4.3294   & 3.0375e-05 &  2.4356  &  1.4551e-05  & 4.8250  \\  \hline
512&   1.7971e-05  & -0.3452 &   1.0095e-05 &  1.5892  &  2.5964e-05 &  -0.8353  \\  \hline
1024&    8.8068e-07  & 4.3509   & 1.8901e-06  &  2.4171 &   9.0832e-07 &  4.8372 \\  \hline
2048 &    1.1262e-06 & -0.3547 & 6.2812e-07  & 1.5894   & 1.6256e-06  & -0.8397  \\  \hline
4096&    5.4989e-08  & 4.3562  &  1.1800e-07 &  2.4123  &  5.6755e-08 &  4.8401  \\  \hline
 &  average &  2.1093 &  average & 2.0778  &average & 2.3974\\ \hline

   \end{tabular}

\end{table}

\begin{example} $\null$
\end{example}
We show example for general self-adjoint two point boundary value problem \eqn{eq1d}:
\eqm
(\beta(x) u_x)_x - \sigma(x) u = f(x).
\enm
The setting is listed  below,
 \eqm
  \beta(x) & = & \left \{ \begin{array}{ll}
 1+ x^2    &   \mbox{if    $0<x< \alf $,} \\ \eqsp
  1.1    &    \mbox{ if $ \alf < x < 1; $}
  \end{array} \right. \\ \eqsp
  f(x) & = & \left \{ \begin{array}{ll}
  2 + 6 x^2   - x^3    &   \mbox{if    $0<x< \alf $,} \\ \eqsp
\dsp     12.2 x^2 -x^5    &    \mbox{ if $ \alf < x < 1, $}
  \end{array} \right. \\
   \sigma(x) &=& x .
  \enm
The solution to the boundary value problem is
\eqm
 u(x)  &=& \left \{ \begin{array}{ll}
 x^2     &   \mbox{if    $0<x< \alf $,} \\ \eqsp
 x^4    &    \mbox{ if $ \alf < x < 1; $}
  \end{array} \right.
\enm

We  discretize  the $ \sigma(x) u(x)$  term as a source term.
In Table~\ref{tab3}, we show the grid refinement analysis for the example. We observe the similar convergence behavior as before, that is, average second order convergence in the infinity norm.

\begin{table}[hptb]
\centering \caption{Grid refinement analysis of the improved harmonic average method for the general self-adjoint interface problem with  $\alf=5/9$.
} \label{tab3}

\vthin

\begin{tabular}{||c||c|c||c| c||c| c||}
\hline
     $N$ &   $\|E \|_{\infty}$ & order     \\  \hline

 32&   3.7309e-04      &       \\  \hline
 64 &  8.0773e-05   & 2.2076 \\  \hline
 128 &   9.8570e-05  & -0.2872  \\  \hline
  256 &  1.6861e-05 &   2.5475 \\  \hline
 512 &  1.7392e-06  & 3.2772 \\  \hline
 1024 &  3.3307e-07 &  2.3845 \\  \hline
  2048 &  8.4962e-08 &  1.9709 \\  \hline
 4096 &  1.8762e-08  &  2.1790 \\  \hline
  & average & 2.0399 \\  \hline
\end{tabular}
\end{table}

\begin{example} A two-dimensional example with a line interface that is parallel to one of axis.
\end{example}
For linear two-dimensional elliptic boundary value problems, the improved harmonic average method works fine if the interface is a line that is parallel to one of axis since the discretization can be done dimension by dimension.
We present numerical results in Table~\ref{tab2d} with the following settings.
 \eqm
  \beta(x,y) & = & \left \{ \begin{array}{ll}
\beta^-   &   \mbox{if    $0<x< \alf $, $0<y< 1$,} \\ \eqsp
  \beta^+    &    \mbox{ if $ \alf < x < 1$,  $<y< 1$;}
  \end{array} \right. \\ \eqsp
  f(x,y) & = & \left \{ \begin{array}{ll}
   -   \beta^-    \left( k_1^2  +1 \right ) \sin(k_1 x)   \cos y   &   \mbox{if    $0<x< \alf $, $0<y< 1,$} \\ \eqsp
\dsp   -   \beta^+  \left( k_2^2 +1 \right )  \cos (k_2 x)   \cos y  &    \mbox{ if $ \alf < x < 1$, $0<y< 1.$}
  \end{array} \right.
  \enm
The analytic solution is
\eqm
  u(x,y) & = & \left \{ \begin{array}{ll}
\sin(k_1 x) \cos y,   &   \mbox{if    $0<x< \alf $, $0<y< 1$,} \\ \eqsp
\cos  (k_2 x) \cos y  &    \mbox{ if $ \alf < x < 1$,  $<y< 1$.}
  \end{array} \right.
 \enm
In the numerical experiments, the Dirichlet boundary condition is used.

\begin{figure}[hpt]
\begin{minipage}[t]{2.7in} (a)\\

\includegraphics[width=0.9\textwidth]{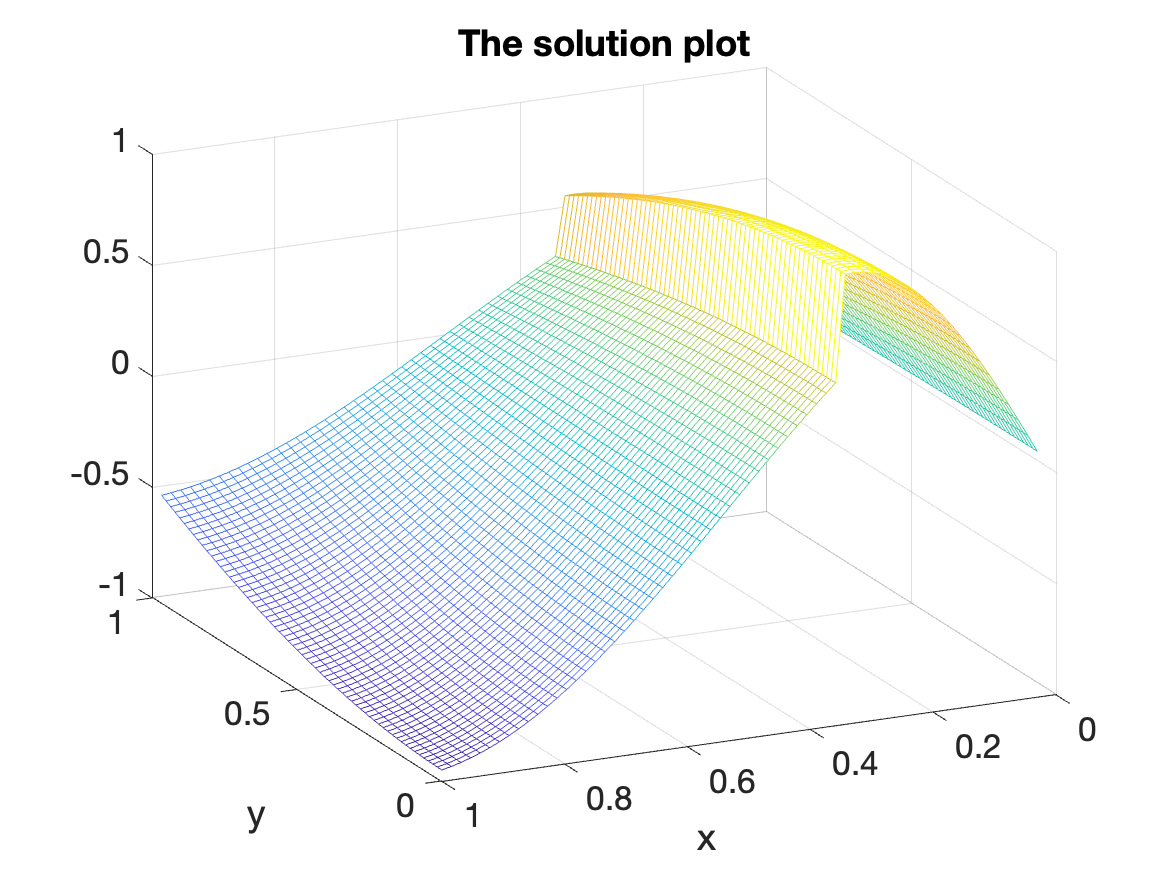}
 \end{minipage}
\begin{minipage}[t]{2.7in} (b)\\

\includegraphics[width=0.9\textwidth]{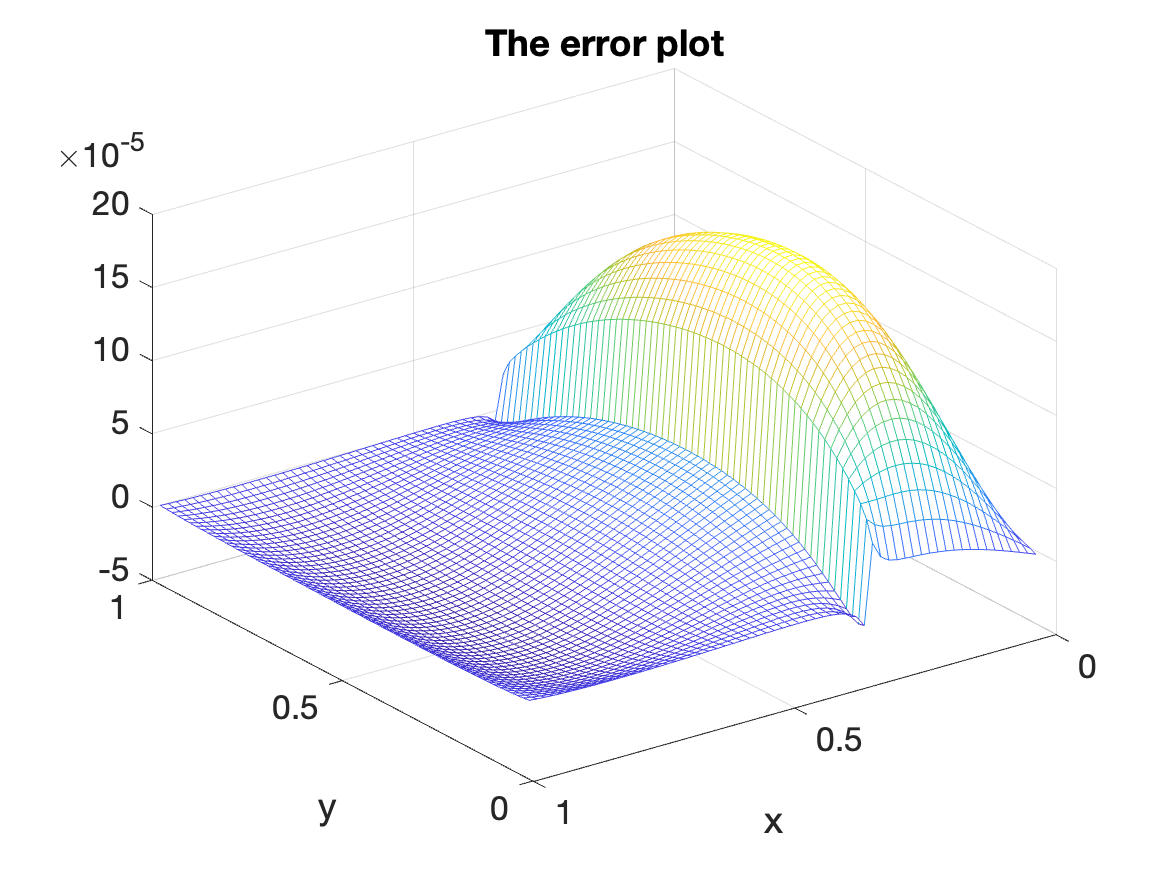}
 \end{minipage}


\caption{(a), the computed  solution when $\beta^-=1.5$, $\beta^+=3$, $k_1=5, k_2=3$ and $n=64 \times 64$.  (b), the corresponding error plot. }\label{fig:2d}
\end{figure}

In Table~\ref{tab2d}, we present a grid refinement analysis for the 2D example with  $\alf=1/3$,  $\beta^-=1.5$, $\beta^+=3$, $k_1=5, k_2=3$. We observe similar behaviors as those in one-dimensional examples. The average convergence order is $2.0207$. If we take the data of $16 \times 16$ away, then the average convergence order  is $2.4952$. In Figure~\ref{fig:2d}~(a), we show the mesh plot of the computed solution with $\alpha = 1/3$, $k_1=5, k_2=3$ using a $64 \times 64$ grid in which the jump discontinuity is captured accurately.  In Figure~\ref{fig:2d}~(b), we plot the errors of the computed solution in which  the plot of error shows where the discontinuity is.

\begin{table}[hptb]
\centering \caption{Grid refinement analysis of the improved harmonic average method for the two dimensional
 interface problem with  $\alf=1/3$,  $\beta^-=1.5$, $\beta^+=3$, $k_1=5, k_2=3$.
} \label{tab2d}

\vthin

\begin{tabular}{||c||c|c||c| c||c| c||}
\hline
     $n\times n$ &   $\|E \|_{\infty}$ & order     \\  \hline
$16\times 16$	& 3.0436e-03 & \\  \hline
$32\times 32$	& 3.8852e-03	&-0.3522 \\  \hline
$64\times 64$	&1.7794e-04	& 4.4485 \\  \hline
$128	\times 128$ & 2.5518e-04 &	-0.5201 \\  \hline
$256	\times 256$ & 1.0953e-05 &	 4.5421 \\  \hline
$512	\times 512$ & 1.6136e-05 &	-0.5589 \\  \hline
$1024\times 1024$	& 6.8190e-07	& 4.5646 \\  \hline
  & average & 2.0207 \\  \hline
\end{tabular}
\end{table}

\section{Conclusions}

This is probably first paper that analyzes the convergence of the harmonic average method for interface problems. The tool is the second-kind Green's function. The analysis shows that the  harmonic average method is second order accurate in the infinity norm when the solution and the flux are continuous (homogeneous jump conditions). Based on the analysis, we proposed the improved harmonic average method that works for more general problems whether the solution and/or the flux have finite jumps or not.
Compared with existing method, the improved harmonic average method is derivative free and has been shown to be second order accurate in the strongest point-wise norm, not the average norm used in the finite element formulation.
Numerical experiments have confirmed the theoretical analysis. Usually, the harmonic average method  has only first order accuracy for general two- or three-dimensional interface problems except the cases that the interface is parallel to one of axis.

\section*{Acknowledgements}
   Zhilin Li was partially supported by a Simons grant 633724.
Kejia Pan is supported by the National Natural Science Foundation of China (No. 42274101).

\section*{Conflict of interest statement}

The authors whose names are listed immediately below certify that they have NO affiliations with or involvement in any organization or entity with any financial interest (such as honoraria; educational grants; participation in speakers’ bureaus; membership, employment, consultancies, stock ownership, or other equity interest; and expert testimony or patent-licensing arrangements), or non-financial interest (such as personal or professional relationships, affiliations, knowledge or beliefs) in the subject matter or materials discussed in this manuscript. 




\bibliographystyle{amsplain}

\bibliography{../../TEX/BIB/bib,../../TEX/BIB/zhilin,../../TEX/BIB/other,harmonic}

%
%

\end{document}

%% file: zmacros.tex
\newcommand{\eqn}[1]{(\ref{#1})}

\newcommand{\eq}{\begin{equation}}
\newcommand{\en}{\end{equation}}
\newcommand{\eqm}{\begin{eqnarray}}
\newcommand{\enm}{\end{eqnarray}}
\newcommand{\eqmno}{\begin{eqnarray*}}
\newcommand{\enmno}{\end{eqnarray*}}
\newcommand{\eql}[1]{\begin{equation}\label{#1}}
\newcommand{\half}{{\frac{1}{2}}}

\newcommand{\goto}{\rightarrow}
\newcommand{\vthin}{\vskip 5pt}
\newcommand{\ignore}[1]{}
\def\bc{\begin{center}}
\def\ec{\end{center}}
\def\bi{\begin{itemize}}
\def\ei{\end{itemize}}
\def\be{\begin{enumerate}}
\def\ee{\end{enumerate}}

\def\alf{\alpha}

\def\reals{{{\rm l} \kern -.15em {\rm R} }}

\def\qquad{\quad\quad}

\def\eqsp{\noalign{\vskip 5pt}}
\def\qquad{\quad\quad}

\def\dsp{\displaystyle}
\newcommand{\eqml}[1]{\eql{#1}\begin{array}{rcl}}
\newcommand{\enml}{\end{array}\end{equation}}